\documentclass[12pt]{amsart}
\usepackage{geometry}                % See geometry.pdf to learn the layout options. There are lots.

\usepackage[parfill]{parskip}    % Activate to begin paragraphs with an empty line rather than an indent
\usepackage{graphicx}
\usepackage{amssymb}
\usepackage{epstopdf}
\usepackage{color}

\newtheorem{theorem}{Theorem}[section]
\newtheorem{prop}[theorem]{Proposition}
\newtheorem{lemma}[theorem]{Lemma}
\newtheorem{cor}[theorem]{Corollary}

\newcommand{\btheorem}{\begin{theorem}}
\newcommand{\etheorem}{\end{theorem}}
\newcommand{\bprop}{\begin{prop}}
\newcommand{\eprop}{\end{prop}}
\newcommand{\blemma}{\begin{lemma}}
\newcommand{\elemma}{\end{lemma}}
\newcommand{\bcor}{\begin{cor}}
\newcommand{\ecor}{\end{cor}}

\theoremstyle{remark}
\newtheorem{remark}[theorem]{Remark}

\newcommand{\bg}{{\mathfrak g }}
\newcommand{\ba}{{\mathfrak a }}
\newcommand{\bn}{{\mathfrak n }}

\newcommand{\bl}{{\mathfrak l }}
\newcommand{\bk}{{\mathfrak k }}
\newcommand{\bq}{{\mathfrak q }}

\newcommand{\bh}{{\mathfrak h }}
\newcommand{\bt}{{\mathfrak t }}

\newcommand{\bp}{{\mathbf p}}
\newcommand{\bc}{{\mathfrak c}}

\newcommand{\bR}{{\mathbb R}}
\newcommand{\bN}{{\mathbb N}}
\newcommand{\bC}{{\mathbb C}}
\newcommand{\bZ}{{\mathbb Z}}

\newcommand{\C}{{\mathbb C}}
\newcommand{\R}{{\mathbb R}}

\title{Branching laws for discrete series of some affine symmetric spaces}
\author{Bent {\O}rsted and Birgit Speh}
\address{B. {\O}rsted, Department of mathematics, Aarhus University, 8000 Aarhus C, Denmark
\newline
Email:
orsted@math.au.dk}
\address{B. Speh, Department of Mathematics, Cornell University, Ithaca 
  NY 14853, USA.  
\newline
Email: 
bes12@cornell.edu}
\thanks{Research by B.Speh partially supported by NSF grant DMS-1500644 }

\begin{document}

\maketitle

\

\begin{center}
{\em Dedicated to B.~Kostant }
\end{center}

\begin{abstract}
In this paper we study branching laws for certain unitary representations. This is done
on the smooth vectors via a version of the {\it period integrals}, studied in number theory, and also
closely connected to the {\it symmetry-breaking operators}, introduced by T.~Kobayashi.
 We exhibit non-vanishing symmetry breaking operators for   the restriction of a 
representation $\Pi$ in the discrete spectrum for real hyperboloids to 
representations of smaller orthogonal groups. In the last part we discuss  
some conjectures for the restriction of representations in Arthur packets 
containing the representation $\Pi$ and the corresponding Arthur-Vogan packets 
to smaller orthogonal groups; these are inspired by  the Gross-Prasad conjectures.  

\end{abstract}

\

\

{Mathematic Subject Classification (2010):} 
Primary 22E46; 
Secondary 
53C30, 
22E30
22E45,
22E50

\section{Introduction}
%{\color{magenta} Start with the number theoretic significance of periods, and special values of L-functions and examples: Rankin Selberg, Ichino-Ikeda , Azai L-%functions, check .}

\

%{\color{red} IMPORTANT Question:  Do we have in our computations anywhere the assumption G = O(p,q) with p even? }

The restriction of a finite-dimensional irreducible representation $\Pi_G$
 of a connected compact Lie group $G$ to a connected Lie subgroup $H$ is a classical problem.
 For example,  the restriction of irreducible representations  of $SU(n+1)$  to $SU(n)$ or of $SO(n+1)$ to 
the subgroup $SO(n)$ 
can be expressed as a combinatorial pattern satisfied by the highest weights of
 the irreducible representation $\Pi_G$ of the large group 
and of the irreducible representations $\pi_H$ appearing in the restriction \cite{We}.
 More generally determining the restriction of an irreducible representation of a compact connected Lie group to a closed connected subgroup
 is transformed to a combinatorial problem  involving highest weights by a famous 
branching theorem by B.~Kostant  \cite{Ko}.  
 
Considerable efforts have been devoted recently to understanding the restriction or branching 
laws for infinite dimensional irreducible 
 representations $\Pi$ of a semisimple
Lie group $G$ to a  symmetric subgroup $G'$ and more recently to other subgroups.
There are many different techniques for dealing with this question. 
The orbit method and geometric quantization considered by B.~Kostant are  
useful tools  to determine the restriction of the representation 
$\Pi$ to a non-compact reductive subgroup $G'$ (see  \cite{K-N}) if the  restriction of 
$\Pi $ is $G'$-admissible, i.e. if it is a 
direct sum of irreducible representations with finite multiplicities
of the subgroup $G'$. For example this is the case if  we restrict  a holomorphic discrete series representation  $\Pi$ to a large subgroup. See  \cite{P}.

\medskip

\bigskip
%In this article 
%our main tool is an integral closely related to  {\it period integrals} of automorphic representations.  
%we do not consider the restriction of discrete series representations of a semisimple group  $G$ i.e irreducible subrepresentations  of $L^2(G)$, 
To determine  explicit branching laws, 
one usually uses specific  models of the representations involved. In the present
paper we work with  representations in the way they arise explicitly
as discrete series of $L^2(G/H)$ of a symmetric space $G/H$, and we refer to those
 as Flensted-Jensen representations. 
We  investigate  the restriction of discrete series representations of symmetric spaces $G/H$, i.e  of irreducible subrepresentations  $\Pi$ of $G$  in  $L^2(G/H)$, to a symmetric subgroup $G'$. % in the "smooth category".
 The restriction of  $\Pi $  to   $G'$ is  not always a direct sum of irreducible representations of 
$G'$ and thus as in \cite{K-S1}  and in \cite{K-S2} for irreducible representations $\pi$ of $G'$
we consider  
\[ \mbox{Hom}_{G'}(\Pi^\infty_{|G'} ,\pi^\infty) .\] 
Here, as the notation suggests, we  work in the smooth 
category and consider  $G'$-equivariant continuous 
homomorphisms  in the topology on the $C^\infty$-vectors introduced by Casselman-Wallach 
\cite{W1} and Bernstein-Kroetz \cite{B-K}. %and we 
%and the multiplicity
%\[ m(\Pi, \pi)=\mbox{dim Hom}_{G'}(\Pi^\infty_{|G'} ,\pi^\infty )\] 
Equivalently  we can consider periods
\[ \mbox{ Hom}_{G'}(\Pi^\infty_{|G'} \otimes (\pi^\vee )^\infty , \bC).\] 
where  $(\pi^\vee )^\infty $ is the smooth part of the contragredient representation of $\pi$.
%but we obtain in special cases criteria for similar representations $W$ of $G'$
%to occur in the restriction $V_{|G'}$.
 
In our case we shall consider a linear functional on $\Pi^\infty \otimes ( \pi^\vee )^\infty$ 
 defined by a period integral over $G'$-orbits on $X=G/H$; namely
$$I = I(v,w) = \int_{X'} v(\dot g' e_p) \cdot w(\dot g' e_p) d\dot g$$where $X' = G'/H', \, H' = G' \cap H$ and $H$ is the stabilizer of the base point $e_p$.  This linear functional is $G'$-invariant provided it converges for all $C^\infty$-vectors. We obtain a branching law by proving that this linear functional is not zero. For discrete
series of a semisimple Lie group a  similar linear functional 
was used  by J. Vargas \cite{Va} to obtain in some cases a branching law 
for unitary restriction of  discrete series representations of $G$ to a subgroup $G'$. We
emphasize that here we work in the category of smooth representations.

% If U and W are smooth representations of G, respectively of G'( in the Casselman Wallach category)  then 
 %  \[       \mbox{Hom}_{G'}( U\boxtimes W^{t}, \bC)  =     \mbox{Hom}_{G'}( U, W) .\].
 
%(see K-S)  where $W^t $ denotes the contragradient representation.
% So a nonzero periods of $V\boxtimes W$  implies that we have a nonzero { Symmetry breaking operators in}
  %             \[  \mbox{Hom}_H( U, W) \]

\bigskip
 In this article we  illustrate the main ideas in one specific example
(making formulas quite explicit). We assume  that 
    \[G= SO_o(p,q) \quad \mbox{with } \quad  p,q \geq 4\ \quad  \mbox{and }
H= SO_o(p,q-1) .\]
and consider representations $\Pi$ in the discrete spectrum of $L^2(G/H)$ of the de Sitter space $X=G/H$ and  representations $\pi^\vee$ in the discrete spectrum of  $L^2(G'/H')$ where 
\begin{enumerate}
\item $G'=SO_o(p-1,q)$ and $ H'= SO_o(p-1,q-1), $ or
\item $G'=SO_o(p,q-1)$ and $H'=SO_o(p,q-2)$
\end{enumerate}
 We discuss  two basic questions for our period integrals,
namely their convergence and their non-vanishing on $\Pi^\infty \otimes (\pi^\vee )^\infty$. 
For this we use the work of M.~Flensted-Jensen \cite{FJ} about discrete series for 
affine symmetric spaces $G/H$, criteria for existence, and
explicit formulas for test functions in the minimal $K$-types. 
We show that for certain pairs of representations our period integrals do not vanish for 
functions  in the minimal K-types and that they define in this case a 
$G'$-invariant linear functional  on the underlying $(\bg,K)$-modules.  Then we show that the period integral 
$I$ in fact defines a continuous $G'$-invariant linear functional on the $C^\infty $-vectors. 
In the last step, i.e.  to show that this linear functional is continuous,  we use that 
the smooth vectors in $\Pi \times \pi$ form an irreducible module for the algebra of 
Schwartz functions on $G, G'$ \cite{B-K}, \cite{W1}. Thus we obtain the "automatic continuity" 
of the period integral $I$.

\medskip

Following the notation in \cite{S1} we shall (see more details below) 
parametrize the representations in the discrete series of  $G/H$, $G'/H'$
 and $G''/H''$ (cases (1) and (2) above)
 by constants $a, \ b, \ c \ \in {\bZ_+}{/2}$ respectively. Note that if $p+q$ is even then $a$ is integral and $b$, $c$
 are half integral. If $p+q$ is odd then $a$ is half integral and $b$, $c$ are integral.

\medskip
\begin{theorem}
 Under the above assumptions
 \begin{enumerate}
 \item
 Suppose that $G'=SO_o(p-1,q)$ and $ H'= SO_o(p-1,q-1).$
The period integral is nontrivial on the minimal K-type of $\Pi_a \otimes \pi_b$ if $b <a$.
\item
Suppose that  $G'=SO_o(p,q-1)$ and $H'=SO_o(p,q-2)$.
The period integral is nontrivial on the minimal K-type of $\Pi_a \otimes \pi_c$ iff $c=a+1/2$.
\end{enumerate}
\end{theorem}

\medskip

%These results are already of separate interest if q=1 and G'=O(p-1,1), It implies that some of the symmetry breaking operators in in Kobayashi-Speh [KS1]. 
%are not trivial on the minimal K-type.

\medskip
 We  can now formulate our branching results. \\
 We say that two finite  sequences of half integers
 \[a_1,\ a_2, \dots ,\ a_m \] 
 \[ b_1,\ b_2,   \dots , \ b_n\]   have  an interlacing property  of finite type if
 \[a_1> b_1 > a_2 > b_2 >  \dots \]
 We say that the sequences have a interlacing property of infinite type 1 if  
  \[b_1 \geq a_1 > a_2 > b_2 >  \dots \]  %\footnote{$b_1$ larger or equal $a_1$ }
 
 Note that the infinitesimal character of the representation $\Pi_a$ is
 \[ (a,0,0,  \dots ,0) +\rho  = (a+\frac{p+q}{2} , \frac{p+q-2}{2},  \dots )\]

 \begin{theorem}
 Under the above assumptions
 \begin{enumerate}
 \item
 Suppose that $G'=SO_o(p-1,q)$ and $ H'= SO_o(p-1,q-1), $
 Suppose that  $\Pi$ and $\pi$ are Flensted Jensen representations with  infinitesimal characters
 \[(a+\frac{p+q}{2} , \frac{p+q-2}{2},  \dots )\]
\[(b+\frac{p+q-1}{2},  \frac{p+q-3}{2},  \dots )\]
If the infinitesimal characters of $\Pi$ and $\pi$ have an interlacing property of finite type then
\[ \mbox{dim Hom}_{G'}(\Pi^\infty_{|G'} ,\pi^\infty )  \not = 0. \] 

  \item Suppose that  $G'=SO_o(p,q-1)$ and $H'=SO_o(p,q-2)$.
  Suppose that  $\Pi$ and $\pi$ are Flensted Jensen representations with the infinitesimal characters
 \[(a+\frac{p+q}{2} , \frac{p+q-2}{2},  \dots )\]
\[(b+\frac{p+q-1}{2},  \frac{p+q-3}{2},  \dots )\]
Then 
\[ \mbox{dim Hom}_{G'}(\Pi^\infty_{|G'} ,\pi^\infty ) \not = 0 \] 
 implies that  $b \geq a+1/2$ and thus
the infinitesimal characters  have an interlacing property of infinite type 1.
\end{enumerate}
 \end{theorem}
 
 \medskip
 \bigskip
 
 In the last section we look at our results from a different perspective using  ideas coming from 
automorphic forms, the  trace formula and the Gross-Prasad conjectures  \cite{GP}, \cite{K-S3}, \cite{N-V}. They suggest  to consider  not  
individual groups and representations but rather packets of representations of $G=SO_o(p,q)$, i.e.
 an Arthur packet  $\mathcal C_{p,q}$ which contains discrete series representation 
$\Pi_s$ of $L^2(G/H)$. Instead of considering one group $SO_o(p,q)$   
it might be helpful to consider a family of groups namely the pure inner forms 
containing $SO_o(p,q)$. Instead of considering a Arthur packet ${\mathcal C_{p,q}}$ of 
representations  it might best to formulate the results in the language of  
Arthur-Vogan packets $\mathcal AV$. For motivation and details see  \cite{GP}, \cite{A-B-V} and section V.
 
 We fix a regular infinitesimal character.
 \[(s,0,\dots ,o) + \rho \]
of a finite dimensional representation of $SO(p+q)$. In \cite{A-J} J. Adams and 
J. Johnson introduced stable packets  of cohomologically induced representations 
which are usually referred to as Arthur packets. The Arthur packet  containing a 
Flensted-Jensen representation  $\Pi_s$ in the discrete spectrum of $SO_o(p,q)/SO_o(p,q-1) $ 
contains also another  representation $\Pi_{as}$ .  The representations are cohomologically
 induced from parabolic subgroups with Levi $SO(2,0) \times SO_o(p-1,q)  $,  
respectively. $SO_o(p,q-1) \times SO(0,2)$ and the representation $\Pi_{as}$ is contained in 
the discrete spectrum of $SO_o(p,q)/SO_o(p-1,q)$.
%(I believe looking at Schlichtkrulls formulas for the Langlands parameter) that the second representation is a anti de sitter rep. Call them Pi_s and \Pi_a.  

We restrict to  2 subgroups $G^1 = SO_o(p-1,q)$ and $G^2 = SO_o(p,q-1)$ which are not in the same 
inner class of orthogonal groups.  For each of these subgroups we have corresponding 
Arthur packets of representations. $L$-functions and  trace formula considerations strongly suggest  
that the
multiplicities depend  on the inner class of the groups $G^i$ and the Arthur packet of representations.
 
Our computations and the results of \cite{K-O}  support this; they show that the restriction of
 exactly one of the representations in our Arthur packet is direct sum, the other one is not admissible.

\medskip
Based on our computations we expect furthermore  a multiplicity 
one result similar to the famous Gross-Prasad conjecture, namely as follows:
 
{\em Suppose  that $\pi$ is a representations of $G^1$ in an Arthur  packet which contains a Flensted-Jensen  representation in the discrete spectrum of $G^1/(H^1 \cap G^1)$. If

                  \[    \mbox{Hom} _{G^1} ( \Pi_s^\infty, \pi^\infty ) \not = 0 \]
then 
                  \[     \mbox{Hom} _{G^1} ( \Pi_{as}^\infty , \pi^\infty ) = 0  \]

and if        \[      \mbox{Hom} _{G^1} ( \Pi_{as}^\infty, \pi^\infty )\not = 0  \]

then        \[         \mbox{Hom} _{G^1} ( \Pi_{s}^\infty, \pi^\infty ) = 0  \]

We  expect that an analogous statement holds for the group  $G^2$.}

\bigskip
{\bf Acknowledgement.}
That this paper is devoted to B. Kostant is only a small token of the gratitude and
debt that we feel and owe for much inspiration and mathematical energy.
B. Speh would like to thank the department of mathematics at Aarhus University for its 
hospitality and support during this research. We also thank T. Kobayashi for some
private communications, and in particular remarks on the topic of branching laws
for the representations considered here.

\bigskip

\noindent

%%%%%%%%%%%%%%%%%%%%%%%%%%%%%%%%%

 \section{Generalities}

In this section we recall the results about  symmetric spaces by M.~Flensted-Jensen \cite{FJ}, Flensted Jenson functions and Flensted Jensen representations \cite{S2}, \cite{S1}. 
We also discuss the restriction of representations and introduce symmetry breaking 
operators as well as periods. The notation  here is inspired by \cite{K-S1}.

%{\color{red}  \bf Caution!!  I have no access to the literature right now , so details and parameter  %maybe wrong, I also have no complete notes and no access to %the pictures of the blackboard }\\

\bigskip
\subsection{Notation and geometric background} Consider the quadratic form 
 \[ Q(X,X)=x_1^2 + ... + x_p^2 -x_{p+1}^2- ....- x_{p+q}^2 \]
and let $G =G_{p,q}=SO_{o,}(p,q) $ be the 
identity component of the automorphism group $O(p,q)$ of the form $Q(X,X)$.

 Recall from the introduction that we assume in this article  $p,q \geq 4$. %as well as $q\geq4$.
 We denote the stabilizer of $e_{p+q}	$ by $H=H_{p,q-1}$, the stabilizer of $e_{1}$ by $G^1=G^1_{p-1,q}$. and the stabilizer of $e_{p+1}$ by $ G^2=G'^2_{p,q-1}  = O(p,q-1)$.
 The maximal compact subgroups of $G,G^1, G^2 $ and $H$ are 
denoted by $K$, $K^1$, $K^2$ and $ K_H$, and we denote the Lie algebras by the 
coresponding small gothic letters. 
 The connected subgroup with Lie algebra $\bR(E_{p+q,p} +E_{p,p+q})$, is denoted by $J$.
 
 \medskip
  The indefinite hyperbolic space $X(p,q)$ is defined by 
 \[ \{\xi \in \bR^{p+q  }| Q(\xi,\xi) =-1 \}  \]
 It is a rank one symmetric space and isomorphic to $G/H$. As a product
 \[ G= KJH \] 
 and 
 \[X(p,q) = (K/M_H) \ J_+\]
where $M_H $ is the intersection of the centralizer of $J$ in $K$ 
with the compact group $K_H $  (see \cite{S1})  and   $J_+ \subset J$ is 
isomorphic to $ \bR^+$. (Strictly speaking, we are here parametrizing an open
dense subset of $X(p,q)$ but this we ignore here and below.)

\medskip
We observe
\begin{lemma}
Under these assumptions \begin{enumerate}
\item  \[ G^1= K^1 J (H\cap G^1)  \quad \mbox{ and } \quad 
   G^2= K^2 J (H\cap G^2) \] 
\item The set  
 \[ \{\xi \in \bR^{p+q  }| Q(\xi,\xi) =1, \xi_p=0 \} = X(p-1,q) \subset X(p,q) \]
 is an orbit  of $G^1$ and isomorphic to $G^1/G^1\cap H$ . The set 
\[ \{\xi \in \bR^{p+q  }| Q(\xi,\xi) =1, \xi_{p+1}=0 \} = X(p,q-1) \subset X(p,q)\]
is an orbit  of $G^2$ and isomorphic to $G^2/G^2\cap H$. 
\item We have 
\[  X(p-1,q) =(K_{G^1}/M_{G^1})J_+\]
respectively
\[  X(p,q-1) =(K_{G^2}/M_{G^2})J_+ \]
  where $M_{G^i}$  i = 1,2 is the intersection of the centralizer of  $J$ in  $K^i$ with the 
compact group $K_{H}=K \cap H$  (see \cite{S1}  8.4.2). 
  \item  $M_G=SO(p-1 )\times  SO(q-1)$ and so \\ $X(p,q) = (S^{p-1}\times  S^{q-1})\bR^+$,   \\
  $X(p-1,q)= (S^{p-2} \times S^{q-1})\bR^+$ and \\ $ X(p,q-1)=(S^{p-1} \times S^{q-2})\bR^+.   $

  \end{enumerate}
\end{lemma}

%\begin{remark} We observe
%\begin{itemize}
%\item The orbit $X(p-1,q)$ is a principal orbit of $G^1$ on $X(p,q)$ .
%\item The orbit $X(p,q-1)$ is not a principal orbit of $G^2$ on $X(p,q)$
%\end{itemize}
%\end{remark} 

\medskip

On $X(p,q) $ we have the  coordinates \[ \Phi (y,y',t)= (y_1 \sinh(t), \dots y_p\sinh(t) , y'_1 \cosh(t),  y'_2\cosh(t),    \dots , y'_q\cosh (t))\]
where $\sum (y_i)^2 = 1$ , $\sum (y_i')^2 =1$ and
\[\int_{X(p,q)} f(x)dx = 
\int_{S^{p-1} \times S^{q-1}} \int_{o}^\infty f(\Phi(y,y',t))\sinh^{p-1} (t  ) \cosh^{q-1 } (t) dt dy  \]

%{\color{red} note to me: check note to Bent }

 \bigskip
\subsection{Analytic background} The universal  enveloping algebra $U(\bg)$ 
acts on the $C^\infty$-functions on the hyperbolic space $X(p,q) $.  The  algebra  
${\mathcal S}(G)$ of Schwartz functions acts by convolutions as well.  (See \cite{B-K} 
and below for the definition of  ${\mathcal S}(G)$, for more details see III.3.)

% For $\lambda \in \mathbb{C}$,
%we consider the eigenspace of the hyperbolic Laplacian
%$\Delta$  on $X(p+1,q)$:
%\[
%\mathcal{S}ol (G/H, \lambda)
%:= \{ f \in \mathcal{C}^\infty (G/H) : \Delta f = \lambda(\lambda-n) f \}.
%\]
%Since $G$ acts isometrically on $G/H$ the solutions
%$\mathcal{S}ol(G/H,\lambda)$ are invariant under $G$ for every
%$\lambda\in\mathbb{C}$. 
For background on the parameters relevant for the discrete series of an
affine symmetric space, see Schlichtkrull \cite{S2}.
%  For $p= 1$ it was first described in \cite{St}.

\medskip

Suppose that  $\lambda $ satisfies the condition 8.8 and 8.9 in \cite{S1} and let 
$F(\lambda ,x)\in  \mathcal{C}^\infty (G/H) $  be the function constructed by M.~Flensted-Jensen (see 7.3 in \cite{S1}). 
In our example 
\[ F(\lambda, \Phi(y,y',t)) = f_a(y') (\cosh (t))^{-\lambda+1 -\frac{p+q}{2}} \]
where $\lambda>0 $ and $a= \lambda-1 +\frac{p-q}{2}$, and $f_a$ is a spherical harmonic
of degree $a$.
Flensted Jensen proved that the Flensted-Jensen  function $F(\lambda ,x)$ has the following properties:

\begin{enumerate}
\item  $F(\lambda ,x)$ is square integrable.
\item  $F(\lambda ,x)$ decays rapidly on $J_+$ in the coordinates $(K/M_H) \ J_+$, i.e. 
it is a Schwartz function as a function on $J_+$ .
\item $F(\lambda, x)$ generates an irreducible submodule $(\bg,K)$-module 
$W(\lambda) \subset \mathcal{C}^\infty (G/H) $.
\item  $F(\lambda ,x)$ as $K$-module  generates the minimal $K$-type.

\item  The unitary completion of the $(\bg,K)$ module $W(\lambda)$ is a unitary representation,
namely it is a unitary irreducible subrepresentation $L^2$--spectrum of $X(p,q)$.
\end{enumerate}

Furthermore  $F(\lambda ,x)$ generates  an irreducible   ${\mathcal S}(G)$ -module ${W}(\lambda)^\infty$ 
which is the Casselman-Wallach realization of the smooth representation of 
$G$ defined by the $(\bg,K)$-module $W(\lambda)$ \cite{B-K}.

\medskip
We refer to the $(\bg,K)$-module $W(\lambda)$ and its completion $W(\lambda)^\infty$ as Flensted-Jensen representations.
We denote the Flensted-Jensen representations of 
$G^1/G^1 \cap H$ by $V(\nu)$ and of $G^2/G^2\cap H$ by $U(\mu)$.

\begin{remark}
By lemma 4.5 in \cite{K2} the functions in the $(\bg,K) $-module of a Flensted-Jensen representation decay at least as rapidly in the  $J_+ $ direction as the Flensted-Jensen function.
\end{remark}

\medskip
\subsection{Restriction of representations} 
 We consider in this article the restriction of a representation 
$\Pi ^\infty$ of a reductive Lie group $G$ to a noncompact reductive subgroup ${G'}$.  
Since the restriction is usually not a   direct sum of irreducible representations  
we consider instead for an irreducible representation symmetry breaking operators in (the smooth category)
\[ \mbox{Hom}_{{G'}}  (\Pi^\infty, \pi ^\infty) \] 
or eguivalently
\[ \mbox{Hom}_{{G'}}  (\Pi _{|{G'}}^\infty\otimes  (\pi^{\vee})^{\infty}), \bC)  \] 
where $\pi ^\vee $ is the contragredient representation of an irreducible representation 
$\pi$ of $G'$.(See  \cite{K-S2}.)

\medskip
Suppose that  $G'$ and $H$ are the fixed points of 2 different involutions of  a reductive group $G$,
both of which commute with the Cartan involution $\theta$. Let $\Pi $ be a representation on the discrete spectrum of $G/H$. Assume that 
\[ \mbox{Hom}_{G'} (\Pi^\infty _{|G'},  \pi^\infty )  \not = 0 .\] 
   
If $G'/H\cap G' \subset G/H $ is not a principal orbit, then the restriction of $\Pi$ may not be a direct 
sum of irreducible representations, and $\pi$ may not be a 
Flensted-Jensen representation in the discrete spectrum of $L^2(G'/G'\cap H)$. 
Consider the following example: Note that for $G= SO_o(p,q) $, $H= SO_o(p,q-1) $ and a 
Flensted-Jensen representation $\Pi$  for G/H we have
\[ \mbox{Hom}_{H}  (\Pi_{|H}^\infty,  \bC) \not = 0 .\] 
For the subgroup $G^2$ which is conjugate to $ H $  in G we also have 
\[ \mbox{Hom}_{G^2}  (\Pi_{|G^2}^{\infty},  \bC)  \not = 0 ,\] but
the trivial representation is not a Flensted-Jensen representation in \linebreak $L^2(G^2/G^2 \cap H)$.

\medskip

\medskip
\subsection{Flensted-Jensen representations of orthogonal groups and their parameters}
We assume in this subsection that $G=SO_o(p,q)$ and $H=SO_0(p,q-1)$ is the stabilizer of $e_{p+q}$.  
We recall the Langlands parameters, the $\theta $-stable parameter and the infinitesimal character  
of the Flensted-Jensen representations  in $L^2(G/H)$ following the exposition  in \cite{S1}.

\medskip
Let $\sigma$ be the involution of $\bg$ with fixed points $\bh$ and let $\bq$ be its -1 eigenspace. 
We have $H \cap K=SO(p) \times SO(q-1)$ and so $K/H \cap K$ is the sphere $SO(q)/SO(q-1)$.  Then $\bR(E_{p+q,,p+1} -E_{p+1,p+q})$ is a maximal abelian subspace $\bt $ of semisimple elements in $\bq \cap \bk$.  

For the rest of the subsection we assume that   $\lambda \in \bt^*$  satisfies the assumptions in 4.1 of \cite{S1}.
  Let 
\[ \mu_\lambda= \lambda+\rho -2\rho_c \in \bt^*\]
where $\rho$ and $\rho _c $ are the half sum of positive roots of $\bt $ on $\bg $ respectively $\bk$. Thus
\[ \mu_\lambda= \lambda+\frac{p-q}{2}-1 \]
is the highest weight vector of a finite dimensional representation $\delta_\lambda$ of K with a 
$(K\cap H)$-fixed vector, i.e.  of a finite dimensional representation of $SO(q)$ with an $SO(q-1)$ fixed vector 
and highest weight $ ( \mu_\lambda,0,\dots,0 )$

\medskip
Let $W(\lambda)$ be the Flensted-Jensen representation with parameter $\lambda $. The lowest $K$-type of $W(\lambda)$ has a minimal $K$-type with highest weight $\mu_\lambda$ \cite{S1} theorem 5.4. The group
                           \[L = G^{\bt} =SO_o(p,q-2) SO(0,2) \]
 is the Levi subgroup of a $\theta$-stable parabolic subgroup and  $\mu_\lambda$ is a 
differential of a one dimensional representation of $L$.  $W(\lambda)$ is obtained by 
cohomological induction from a $\theta$-stable parabolic with Levi $L$ and the character $\lambda$. 
For more details see IV.1.  Following the terminology of \cite{K-S2}  we refer to  $(L,\lambda) $ as the $\theta$--stable parameter of the Flensted-Jensen representation  $W(\lambda)$ and its completion $W(\lambda)^\infty$. Thus we conclude
 
 \begin{prop}
 Let $\bc$ be a Cartan subalgebra of $\bg \otimes_{\R} \C$ .
We choose the usual set $\Delta$ of roots and  positive roots $ \Delta^+$.
 Then the infinitesimal character of $W(\lambda)$ is 
 \[ (\lambda, 0,0 \dots, 0) +\rho \]
 
   \end{prop}

 \medskip
 
 Let $\ba$ be a maximal abelian subspace in $\bl \cap \bp $ and $A$ a connected 
subgroup with Lie algebra $\ba$. So $A$ has dimension 
$min(p,q-2)$. $P_L = M_L A N_L $ is the cuspidal parabolic subgroup of $G$ with $A$ as its split component.  If $\lambda $ is positive and large enough, then 
 $W(\lambda)^\infty $ is the Langlands subquotient of a principal series representation 
$I(P_L,\delta,s)$ induced from $P_L$ with minimal K-type $\delta_\lambda$, see \cite{S2} 
Theorem 6.1, i.e. it has the Langlands parameter $P_L,\delta,s$. We can now prove
 
 \medskip
 \begin{lemma} 
The representation $W(\lambda)$ is self-dual, i.e. it is isomorphic to its contragradient  
$(W(\lambda)^\vee$.
\end{lemma}
{\em Proof.} 
 The map $\Pi \rightarrow \Pi^\vee$  is an involution on the set of irreducible representations  of a 
semisimple group $G$ which was investigated in \cite{A} and \cite{A-V}. The infinitesimal 
character of $\Pi^\vee$ is a Weyl group conjugate of the negative of the infinitesimal character of $ \Pi$.
 Thus  irreducible  representations and their contragradients have the same 
infinitesimal character  if the  Weyl group contains $-1$ \cite{A-V}.  This is satisfied if $p+q$ is odd.
 If $p+q$ is even then 
 \[ (-\lambda, 0,0 \dots, 0) - \rho \]
 is conjugate to 
  \[ (\lambda, 0,0 \dots, 0) +\rho \]
and thus $W(\lambda)$ and $W(\lambda)^\vee$ have the same  infinitesimal character. Since the minimal $K$-type 
is self-dual both representations also have  the same minimal $K$-type. By 
Theorem 5.4 in \cite{S2}
  the minimal $K$-type and  the infinitesimal character determine uniquely the 
Langlands parameter of an irreducible representation.  
\qed.

\

 \section{ Period integrals.} 
In   this section we  assume  that $G=SO_o(p,q)$.
  We consider the restriction of a Flensted-Jensen representation 
${W}(\lambda)^\infty$  to $G^i$, $i = 1, 2$ and find linear nonzero functionals in 
 \[ \mbox{Hom}_{G'}(W(\lambda)^\infty_{|Gi}\otimes (\pi^\vee)^\infty, \bC) \]
where $\pi^\infty =U(\mu)^\infty $ or $ \pi^\infty ={V}(\nu)^\infty$ is a 
Flensted-Jensen  representation generated by a Flensted-Jensen function 
$F(\nu,x)\in C^\infty (G^1/G^1\cap H )$, respectively $  F(\mu,x)\in C^\infty (G^2/G^2\cap H )$.  
For this we consider period integrals of test functions in the lowest 
$K$-type of ${W}(\lambda)^\infty \otimes (\pi^\vee )^\infty$
 and show in III.2 that the period integrals defined in the introduction 
and in II.3 converge   for some pairs of representations  and are nonzero  on the underlying 
$(\bg,K)$-module.  In III.3 we show that the period integrals also converge on the 
$C^\infty$--vectors. In III.4 we show that  
it defines a continuous  linear functional on the Fr{\'e}chet space using the results of Bernstein-Kroetz.

 \medskip
 \subsection{Branching and period integrals for compact orthogonal groups}
 Consider a compact group $K$,  symmetric subgroups $\ F, \  K', \ F' $ so that 
 
  \[
 \begin{array}{ccc}
K' &  \subset &  K  \\
\cup &          & \cup \\
F'=K' \cap F& \subset & F
 \end{array}.  \]
 
  and  the manifolds $Y= K/F$ and $Y' =K'/F'$.
  
 Let $\phi $ be a $K$-finite function in $C^\infty(Y)$ and  $\phi'$ a $K'$-finite function in $C^\infty (Y')$. We consider the period integral 
 \[
 I_{K \cap F} = \int_{Y'} \phi \cdot \phi'
 \]
and we  want to understand for which pairs of representations it is nonzero in our case.

We assume now
\begin{enumerate}
\item  $K=SO(q ),  F=SO(q-1)$ the stabilizer of the last coordinate. 
\item $K'=SO(q-1), F'=K'\cap  F = SO(q-1)$
\end{enumerate} 
 
The manifolds are $Y = K/F,$ $Y' =K'/F'$ are a sphere and the equator.  Recall that $L^2(X) $ is a direct sum of irreducible representations $\pi_a$ of $SO(q)$ with highest weight
(a,0,\dots, 0). We obtain an embedding of $ \pi_a$ into $L^2(X) $ by 
considering the coefficients $ f_a(k)=(v, \pi_a(k)v_o) $ where $v_o$ is the 
$F$-fixed vector. The classical branching rules for  
the restriction of $\pi_a$ to $SO(q-1)$ show that 

                                                                   \[  (\pi_a)_{|SO(q-1)} = \bigoplus _{0\leq   b \leq a} \pi_b. \]

 $\pi_a$ is also an irreducible  $U(\bk) $-module.     This implies
 
 \begin{prop}
 Let $v_a$ be the highest weight vector of $\pi_a$and $v_b$ a highest weight vector of 
$\pi_b$ with $0 \leq b \leq a$. There exists a 
$D_{a,b}\in U(\bk)$ so that 
                                 \[ \int_{Y'} D_{a,b}f_a(k' )\cdot f_b(k')d \dot k'\not =0. \]
 
 \end{prop}       
 {\em Proof:}  We may assume that $D_{a,b} f_{a} $ is a multiple of $f_{b}$.  \qed

 \subsection{Periods for the $(\bg,K)$--module $W(\lambda) \otimes U(\mu)$}

 Let $\lambda, \mu$ be the parameter of the Flensted-Jensen 
representations of   $G$ and $G^2 $ and  let $F(\lambda ,x)$, $F(\mu ,x')$  the corresponding 
 Flensted-Jensen functions. We consider the period integrals
\[ {\mathcal PI}_{p,q-1}(\mu,\lambda):= \int_{X(p,q-1)} F(\lambda ,y) F(\mu ,y) dy \]

or more generally 

\[ {\mathcal PI}^{D_K}_{p,q-1}(\mu,\lambda):= \int_{X(p,q-1)} D_K F(\lambda ,y) F(\mu ,y) dy \]
for $D_k \in U(\bk)$.
Given $\lambda, \mu $ we consider the following two problems:
\begin{enumerate}
\item  For which $\mu $  exists $D_k \in U(\bk)$ so that the integral 
$ {\mathcal PI}_{p,q}(\mu,\lambda)$ is convergent and non-zero. 
\item Does this period integral define a $G^i$-invariant linear functional of $W(\lambda) \otimes V(\mu)$?
\end{enumerate}
 
 \medskip
 
Given the  assymptotics of the measure and of $F(\lambda)$ and $F(\mu)$ we get convergence if
\[-\lambda +1 - \frac{p+q}{2} -\mu +1- \frac{p+q-1}{2} +p-1+q-2 < 0 \]
and so
$\lambda +\mu> -\frac{1}{2} .$
To obtain nonvanishing of the period integral on the minimal $K$-type we need by III.1 
\[ 0 \leq \mu-1 + \frac{p-q+1}{2} \leq \lambda -1+ \frac{ p-q}{2}. \]
Thus we have 
\[
\mu +\frac{1}{2} \leq \lambda.\]

Since the radial part of the integral over $X(p,q-1)$ only involves
a power of $cosh(t)$, this proves (1) in the 

\medskip

\begin{prop}  

  \begin{enumerate}
\item  Suppose that $G'= G^2= SO_o(p.q-1)$ and that $W(\lambda)$ and $U(\mu) $ are the $(\bg,K)$-modules of Flensted-Jensen representations  of G/H and $G^2/H\cap G^2$.

The period integral converges if $\lambda +\mu >-\frac{1}{2}$. 
It doesn't vanish on the minimal $(K \otimes K')$-type of the $(\bg \otimes \bg', K \otimes K')$ modules $W(\lambda)\times U(\mu)$ if \[
\mu +\frac{1}{2} \leq \lambda . \]

\item Suppose that 
$G'=G^1=SO_o(p-1,q) $ and that $W(\lambda)$ and $V(\nu) $ are the $(\bg,K)$-modules of Flensted-Jensen representations  of G/H and  of  $G^1/(H\cap G^1)$. The period integral converges if $\lambda +\nu >-\frac{1}{2}$ and  it doesn't vanish on the minimal K-type
if \[ \nu = \lambda+1/2 .\]
\end{enumerate}
\end{prop}
{\em Proof:} 
In (2) we  consider the period integral
\[ \int_{X(p-1,q)} D_K F(\lambda ,y) F(\nu ,y) dy \] and proceed as in (1).

\qed

\medskip

 \begin{cor}
  \begin{enumerate}
\item The $(\bg,K) $--module $W(\lambda) \otimes V(\nu)$  has a nontrivial $G^1=SO_o(p-1,q)$--invariant linear functional if
\[ \mu = \lambda+1/2 .\]

\item  The $(\bg,K) $--module $W(\lambda) \otimes U(\mu)$  has a nontrivial $G^2=SO_o(p,q-1)$--invariant linear functional if
\[
\mu +\frac{1}{2} \leq \lambda.\]

\end{enumerate}
\end{cor}
 {\em Proof:}  By  II.3 all functions in $W(\lambda)$ respectively in $V(\nu)$ and $\ U(\mu)$ decay at least as fast as the Flensted Jensen functions. Hence the integral converges for all functions in  
 $W(\lambda ) \otimes V(\nu)$, respectively $W(\lambda) \otimes U(\mu)$.
 \qed

\medskip
\subsection{Periods on ${W}(\lambda)^\infty \otimes {V}(\mu)^\infty $}

\subsubsection{ A linear functional on the vector space  ${W}(\lambda)^\infty \otimes {V}(\mu)^\infty$}

We assume again that G=$SO_o(p,q) $ and $H= SO_o(p,q-1)$ and we first ignore the topology on ${W}(\lambda)^\infty \otimes {V}(\mu)^\infty$.

 We recall from \cite{B-K} and \cite{W1} the definition of the algebra ${\mathcal S}(G)$ of Schwartz functions: 
Let $\sigma$ be a finite dimensional faithful representation with a 
$K$-invariant inner product. 
Define for $g \in G$ a scale
\[ s(g) =||g|| = \mbox{tr}(\sigma(g) \sigma(g) ^\vee) +\mbox{tr}( \sigma (g^{-1}) \sigma (g^{-1})^\vee) \] 
and  we denote by $L$ and $R$ the left and right regular representation.
We call the space of smooth functions (decaying  in all derivatives)
\[ \{ f \ | \ \mbox{sup}_{g \in G}||g||^r  |L(X)R(Y)f(g) |  \ < \infty \} \]
for all $X,Y \in U(\bg)$, $r\in \bN$
the Schwartz algebra ${\mathcal S}(G)$, and give it the topology
via the corresponding seminorms. By definition
$s(g) =s(g^{-1})$, $s(g_1 g_2) \leq s(g_1)s(g_2)$ and  similarly
$s(g_1 g_2^{-1})^{-r} \leq s(g_1)^{-r} s(g_2)^r$. We deduce that if 
$f \in {\mathcal S}(G)$ then 
the function $f(g^{-1})$ is also in ${\mathcal S}(G) $. Furthermore 
for $f \in {\mathcal S}(G)$ and $w(g) \in W(\lambda)$  the integral 
\[I_f (x) = \int_G f(g) w(g^{-1}x) dg \]
converges and is in $C^\infty (X(p,q))$. The map 
\[T_w: f  \rightarrow I_f\]
 defines a homomorphism of the Fr{\'e}chet  convolution algebra of 
Schwartz functions ${\mathcal S}(G)$  onto     ${W}(\lambda)^\infty$ 
 \cite{B-K}. Here we use deep results from the theory of globalizations
of Harish-Chandra modules. So we have the identification, also as topological spaces
via the open mapping Theorem,
 \[ W(\lambda)^\infty  = {\mathcal S}(G)/ \mbox{ker }T_w \]  Furthermore  
${W}(\lambda)^\infty$ is an irreducible ${\mathcal S}(G)$-module. Hence to estimate the 
growth of a function in ${W}(\lambda)^\infty$  
we have to estimate the growth of functions   of the form 
$\int_G f(g) F(\lambda,g^{-1}x) dg $ or equivalently $\int_G f(g) F(\lambda, g x) dg $ with $f \in {\mathcal S}(G)$.
For this we make the following explicit estimates which they seem to depend on having classical
groups and it is not clear whether similar estimates can be made in say exceptional groups.
\medskip
For $g \in G $ we write in block form
 \[ g = \left( \begin{array}{cc}
                 A_1 & A_2\\
                 A_3 & A_4
                 \end{array}  \right)  \]
where \[ A_4A_4^t-A_3A_3^t = I \] and 2 more relations. In the coordinates 
\[ (x',x'') = (y'sinh~ t,y''cosh~t)\] where $|y'| = |y''| =1$ (Euclidian norms) on $X(p,q)$ (see II.1) we have
\[ |A_3x'+A_2x''|^2 - |A_3x' + A_4x''|^2 = -1.\]
Now
$ |A_3x' +A_4x''|$ never vanishes  and there exists a constant $C$ so that 
\[F(\lambda,gx)\ \leq C|A_3x' +A_4x''|^{-a}\] (see II.2) so that we need to understand the asymptotics of 
\[|A_3x' +A_4x''| =
  (1/2) e^t |A_3(1-e^{-2t})y' +A_4(1+e^{2t})y''| .\]
  Thus
  we  estimate from below
  \[Q= |A_3y' + A_4 y''| \]
 with $|y'| = 1 = |y''|.$ We note that the sets $\{A_4y''| \ |y''|^2 =1 \} $ and $\{ A_3y'| \ |y'| =1 \}$ are ellipsoids that 
never intersect. We can find their axes by looking at the quadratic forms 
 \[ |A_4y''|^2 = (A_4A_4^{t} y'',y'' ) \quad \mbox{and} \quad |A_3y'|^2 =(A^{t}_3 A_3 y',y') \]
with eigenvalues  
\[ \lambda_1 \geq \lambda_2\dots\geq \lambda_q \geq 1,\] respectively
\[\mu_1 \geq   \mu_2 \dots \geq 0.\]
Note that $A_4^tA_4$ and $A_4A_4^t $ have the same non-zero eigenvalues and 
similarly $A^t_3A_3$ and $A_3A_3^t$. Hence $\mu_1 = \lambda_1-1, \mu_2 = \lambda_2-1$ etc. Note that for  j sufficiently large j
some  eigenvalues $\mu_j$ may be  zero.

\begin{lemma}
 The 2 ellipsoids have parallel axes.
 \end{lemma}
 \proof
 Note that for $\lambda_1$ and its eigenvector $ A_4^tA_4 e_1 = \lambda_1 e_1$ implies that $|A_4e_1|^2 = \lambda_1 $ and thus $A_4e_1$ is the largest axis of the large ellipsoid. We also have
 \[ A_4 A_4^t A_4   e_1= \lambda _1 e_1 \]
 so that 
 \[A_3A_3^tA_4 e_1 = (\lambda_1-1)A_4e_1 \]
 comparing  to the similar
 \begin{eqnarray*}
 A_3^tA_3 v_1       &  =  & \mu_1 v_1\\
 A_3A_3^tA_3 v_1 &  =  & \mu_1 A_3 v_1
 \end{eqnarray*}
 we deduce that that the axis $A_3 v_1$, respectively $A_4 e_1$, of the smaller and the larger 
ellipsoid are proportional if the multiplicity of the eigenvalue is one.  If the multiplicity  
is higher than one, then can identify the eigenspaces and the axes of the
ellipsoids this way. For the next eigenvalue $\lambda_2 = \mu_2 + 1$ the axis
$A_4e_2$ is perpendicular to $A_4e_1$ and similarly $A_3v_2$ is perpendicular to
$A_3v_1$, and again $A_4e_2$ and $A_3v_2$ are linearly dependent.
 \qed 
 
 So can conclude
 
 \begin{lemma}
 Under the above assumptions   
 \[Q= |A_3y' +A_4y"| \geq (1/4) (\mbox {tr } A_4^tA_4)^{-1}\]
 \end{lemma}
 \proof Note
\[ Q = |A_3y' +A_4y"|      \geq \mbox{min}_i (\sqrt {\lambda_i} -\sqrt{\lambda_i-1})\]
So in particular $Q \geq \mbox{min}_i \frac{1}{2 \sqrt{\lambda_i}}$ and thus
\[Q^2 \geq(\lambda_1 +\lambda_2 + \dots)^{-1}.\]
\qed

This gives, using the scales defining ${\mathcal S}(G)$ the estimate for a typical 
function \[ F(x) = \int_{G/ H } f(g)F(\lambda,gx) dg \in V^\infty \] from above

\begin{prop}
Every $F  \in  V(\lambda)^\infty = {\mathcal S(G)}F(\lambda)$ satisfies 
\[ |F(x) | \leq Ce^{-at}\]
The constant C depends on $f \in {\mathcal S}(G) $ and $C \geq \int_G |f(g)| s(g)^a da $.
\end{prop}
\proof 
For the defining representation  the Hilbert-Schmidt norm of
$g$ dominates the Hilbert-Schmidt norm of $A_4$. Hence the scale $s(g)$ dominates the Hilbert-Schmidt norm $\mbox{tr  }A^t_4 A_4$
and asymptotically for $t\rightarrow \infty$ 
\[| A_3(1-e^{-2t})y' +A_4(1+e{-2t})y''|\] approaches $|A_3y'+A_4y''| $.
\qed

\medskip

In a similar way we can estimate   the growth of functions in $V(\nu)^\infty$ and $U(\mu)^\infty$.
 Hence we can conclude

\medskip
\begin{cor}
 The period integral 
 \[I(f,f')= \int_{G^i/H^i} f\cdot {f' }dg' \]
 of functions $f,\ f'$ in $W(\lambda)^\infty$, $V(\mu)^\infty$ respectively $U(\nu)^\infty$
converges if the integral $I$ over the corresponding Flensted-Jensen functions converges and it 
defines a linear functional $I$ on $W(\lambda)^\infty \otimes V(\mu)^\infty$, and 
$W(\lambda)^\infty\otimes U(\nu)^\infty$ respectively. 
\end{cor}

\medskip
\subsubsection{Continuity of  the period $I$}
To complete the proof it remains to show that the  linear functional  
\[{\mathcal I}:W(\lambda)^\infty \times V(\mu)^\infty \rightarrow \bC\] is continuous.

\medskip
Note the map \[{\mathcal S}(G) \rightarrow V(\lambda)^\infty \]
\[f \rightarrow  I_f(x) = \int_G f(g) F(\lambda,g^{-1}x)dg \]
is continuous. So if $f_n \rightarrow f$ in ${\mathcal S}(G)$ then the estimates imply that 
$I_{f_n} \rightarrow I_f$ pointwise and dominated by $C \, e^{-at}$ for some $C$.

Hence for fixed $F' \in V(\mu)^\infty$
\[I_n=  I( \int_G f_n(g) F(\lambda,g^{-1}x)dg,F') \rightarrow  I( \int_G f(g) F(\lambda,g^{-1}x )dg, F')\]
for $n \rightarrow \infty$. Hence we have continuity in the first variable. A similar argument 
shows that we have continuity in the second variable. Since we are working in 
Fr{\'e}chet spaces the period integral is continuous in both variables.

\medskip
Thus we conclude

\medskip
\begin{theorem}
\begin{enumerate}
\item  The representation ${W}(\lambda)^\infty \otimes {V}(\nu)^\infty$
 has a nontrivial $G^1=SO_o(p-1,q)$--invariant linear functional which is continuous in the Casselman-Wallach topology if 
 \[ \mu = \lambda+1/2. \]
\item The representation 

  ${W}(\lambda)^\infty \otimes {U}(\mu)^\infty$
 has a nontrivial $G^2=SO_o(p,q-1)$--invariant linear functional which is continuous in the Casselman-Wallach topology if
\[
\mu +\frac{1}{2} \leq \lambda\]
\end{enumerate}
\end{theorem}

\medskip
\begin{remark}
This result shows that the linear functional defined by a period integral on the $(\bg,K)$-module 
$W(\lambda)\boxtimes V(\nu)$ extends "automatically" to a continuous linear 
functional on the completion in the spirit  of the automatic continuity theorems by Casselman \cite{C} and Delorme, van den Ban \cite{DvB}.

\end{remark}

\section{Branching of ${W}(\lambda)^\infty$ to $G^1$ and $G^2$.}

 In this section we restate the results of the previous section using the terminology of symmetry 
breaking in terms of Langlands parameters as well as the interlacing property of the infinitesimal characters.

  \subsection{Review of cohomological induced representations }
  We summarize 
 cohomological parabolic induction.  
A basic reference is the book by A.W. Knapp and D.~Vogan \cite{K-V}.  
We begin with a {\it{connected}} real reductive Lie group $G$.  
Let $K$ be a maximal compact subgroup, 
 and $\theta$ the corresponding Cartan involution.  Let $H_K$ be a Cartan subgroup of K.
Given an element $X \in \mathfrak{h}_K$, 
 the complexified Lie algebra
 ${\mathfrak{g}}_{\mathbb{C}}={\operatorname{Lie}}(G) \otimes_{\mathbb{R}}
{\mathbb{C}}$ is decomposed into the eigenspaces
 of $\sqrt{-1}\operatorname{ad}(X)$, 
and we write 
\[
   {\mathfrak{g}}_{\mathbb{C}}
   ={\mathfrak{u}}_- + {\mathfrak{l}}_{\mathbb{C}} + {\mathfrak{u}}
\]
 for the sum of the eigenspaces 
 with negative, zero, 
 and positive eigenvalues.  
Then ${\mathfrak{q}}:={\mathfrak{l}}_{\mathbb{C}}+{\mathfrak{u}}$
 is a $\theta$-stable parabolic subalgebra 
 with Levi subgroup 
\begin{equation}
\label{eqn:LeviLq}
   L =\{g \in G: {\operatorname{Ad}}(g) {\mathfrak{q}}={\mathfrak{q}}\}.  
\end{equation}
The homogeneous space $G/L$ is endowed
 with a $G$-invariant complex manifold structure 
 with holomorphic cotangent bundle $G \times_L {\mathfrak{u}}$.  
As an algebraic analogue of Dolbeault cohomology groups
 for $G$-equivariant holomorphic vector bundle over $G/L$, 
 Zuckerman introduced a cohomological parabolic induction functor
 ${\mathcal{R}}_{\mathfrak{q}}^j(\cdot \otimes {\mathbb{C}}
{\rho({\mathfrak{u}})})$ ($j \in {\mathbb{N}}$) from the category
 of $({\mathfrak{l}}, L \cap K)$-modules
 to the category of $({\mathfrak{g}}, K)$-modules. 
 
 \medskip
 Now assume that  $G=SO_o(p,q)$ and
                           \[L = G^{\bt} =SO_o(p,q-2) SO(0,2) \]
 is the Levi subgroup of a $\theta$-stable parabolic subgroup $\bq= \bl + \bn$ and  $\lambda$ is a differential of a one dimensional representation $\lambda$ of $L$ . We assume that the positive roots have positive inner product with the roots in $\lambda$. Then $ W(\lambda) $ is cohomologically induced from ($L, \lambda$). Since we assume 
that $p$ and $q$ are larger than $4$, $L$ is not compact and thus the representations $W(\lambda)$, $V(\nu)$ and $U(\mu)$ are not tempered.

  \subsection{Branching Theorem}
  
 \begin{lemma}
 \begin{enumerate}
\item The restriction of ${W}(\lambda)^\infty $ to $G^1$ is admissible.
 
\item The restriction of ${W}(\lambda)^\infty$ to $G^2$ is not admissible.
\end{enumerate}
 \end{lemma} 
 \proof It follows from the tables in \cite{K-O}.  \qed

 \medskip
 
 Recall that in the Casselman-Wallach realization of irreducible representations 
$\Pi$, $\pi$ of $G$, respectively of $G'\subset G$
  \[ \mbox{Hom}_{{G'}}  (\Pi^\infty_{|G'}, \pi ^\infty) \] 
and 
\[ \mbox{Hom}_{\hat{G'}} ( \Pi _{|{G'}}^\infty \otimes ( \pi ^\vee)^\infty, \bC)  \] 
are isomorphic,
where $\pi ^\vee $ is the contragredient representation of an irreducible representation 
$\pi$ of $G'$. So using Lemma II.4 we can summarize the results of the last 
section in the  \underline{\em Branching Theorem}

\medskip

\begin{theorem}
 Under the previous assumptions 
 \begin{enumerate}
 \item If  $G'= G^2= SO_o(p.q-1)$ and 
\[
\mu +\frac{1}{2} \leq \lambda\] then 

  \[\mbox{Hom}_{G^2}({W}(\lambda)^\infty _{|G^2 } , {V}(\mu)^\infty )\not = 0\]

\item
If $G^1=SO_o(p-1,q) $ and 
\[ \nu = \lambda+1/2 \] then
  \[\mbox{Hom}_{G^1}({W}(\lambda)^\infty _{|G^1} , {V}(\nu)^\infty )\not = 0\]

\end{enumerate}
\end{theorem}

\medskip

\begin{remark}
Using different techniques similar results were obtained in \cite{L}.
\end{remark}

\bigskip

\subsection{Different formulations of the Branching Theorem}
\subsubsection{A graphic formulation   of   the Branching Theorem.}
 
 The restriction of the representation $W(\lambda)$ to $G^2$ can be represented by the following diagram. The representations of G are in  the first row, the representations of $G^i$ in the second row.
 Representations are connected by an arrow if their multiplicity is one.
 
\underline{ Restriction to $G^2$:} Let that $V(\mu)$ be a representation of $G^2$ with $\mu +\frac{1}{2} \leq \lambda$. Then
 \begin{eqnarray*}
  id  & \ & {W}(\lambda)^\infty \\
 \downarrow & \swarrow   & \downarrow\\
 id & \  & {V}(\mu)^\infty
  \end{eqnarray*}
 
\underline{ Restriction to $G^1$}.
Suppose that  $\nu =\lambda+1/2$ 
 \begin{eqnarray*}
  id  & \ & {W}(\lambda)^\infty \\
 \downarrow &   & \downarrow\\
 id & \  & {V}(\nu)^\infty
  \end{eqnarray*}
  \proof
 In the first case the arrow  in the restriction to the trivial representation is the $H=SO_o(p,q-1)$ linear functional defining by the embedding of ${W}(\lambda)^\infty$ into the Schwartz space on $G/H$.
 \qed

{\bf Remark:}  The restriction of ${W}(\lambda)^\infty$ follows the same pattern as the restriction of representations of height $1$ for  $O(n,1)$ to $O(n-1,1)$ \cite{K-S4}.

\medskip

\subsubsection{The Branching Theorem using interlacing patterns }

Interlacing patterns of highest weights of finite dimensional representations and of  $\theta$-stable parameters of irreducible representations of O(n,1) have been used successfully to express the restriction of irreducible representations as well as in the Gross Prasad conjectures. \cite{K-S3}, \cite{K-S4}.

 We say that two finite  sequences
 \[a_1,\ a_2, \dots ,\ a_m \] 
 \[ b_1,\ b_2,   \dots , \ b_n\]   $n = m-1,$ or $n=m$ have  an interlacing property  of finite type if
 \[a_1> b_1 > a_2 > b_2 >  \dots \]
 We say that the sequences have a interlacing property of infinite type 1 if  
  \[b_1> a_1 > a_2 > b_2 >  \dots \]
 
 Note that the infinitesimal character of the representation ${W}(\lambda)^\infty$ is
 \[ ( \lambda,0,0,  \dots ,0) +\rho  = (\lambda+\frac{p+q}{2} , \frac{p+q-2}{2},  \dots )\]
Here we choose a root system so that and positive with respect to the positive root with respect to 
$( \lambda,0,0,  \dots ,0)$ and the roots of $so(q) \oplus so(p)$.

 \begin{theorem}
 Under the above assumptions
 \begin{enumerate}
  \item 
 Let  $\Pi={W}(\lambda)$ and $\pi={V}(\nu)$ be Flensted Jensen representations of G/H respecitively $G^1/(H \cap G^1)$ with the infinitesimal characters
 \[(\lambda+\frac{p+q}{2} , \frac{p+q-2}{2},  \dots )\]
\[(\nu+\frac{p+q-1}{2},  \frac{p+q-3}{2},  \dots )\]
Then 
\[ \mbox{dim Hom}_{G^1}(\Pi_{|G^1}^\infty ,\pi ^\infty) \not = 0 \] 
 implies that $\nu> \lambda $ and thus
the infinitesimal characters  have an interlacing property of infinite type 1.

 \item
 Suppose 
 that  $\Pi ={W}(\lambda)$ and $\pi={U}(\mu )$ are Flensted Jensen representations of $G/H$ respectively $G^2/(G^2 \cap H)$ with the infinitesimal characters
 \[(\lambda+\frac{p+q}{2} , \frac{p+q-2}{2},  \dots )\]
\[(\mu+\frac{p+q-1}{2},  \frac{p+q-3}{2},  \dots ).\]
If the infinitesimal characters have an interlacing property of finite type,
then 
\[ \mbox{dim Hom}_{G^2}(\Pi_{|G^2}^\infty ,\pi ^\infty) \not = 0 .\] 

 \medskip

\end{enumerate}
 \end{theorem}

%Furthermore by private communication by Y. Oshima

%\begin{lemma}
% Suppose that  $G'=G^1=SO_o(p-1,q) $ and $W(\lambda) $ is of discrete series type. Then
%\[\mbox{Hom}_{G'}(\hat{W}(\lambda) , \hat{V}(\mu) )= 0\]
 %if  $\mu <  \lambda+1/2 $.
% \end{lemma}
%{\em A sketch of a proof}

%$\lambda$ is the highest weight of the minimal K-type in $W(\lambda)$ respectively $\hat{W}(\lambda )$.
%\qed

\medskip

% We now recall the description of the representations in the discrete spectrum of $L^2(G/H)$ constructed by Flensted Jensen, see [S].
% Fix a fundamental Cartan subgroup $C$ and 
% a set of  positive roots $\Sigma^+$ of $\bc_\bC, \bg_\bC$ where all simple roots but one are compact.
 % We recall the $\theta$-stable parameter of the representations in the discrete spectum of $X(p+1,q)$ from [K]. 
%All the representations  in the discrete spectrum of $L^2(G/H)$ are cohomologically induced from a $\theta$- stable parabolic subalgebra $\bq_\bC $.
% The  parabolic subalgebra $\bq_\bC$ is a maximal parabolic with a Levi subgroup
%  \[ L= SO_o(p-2,q-2  ) SO(2) \]
% The nilradical $\bu $ of $\bq$ is defined by the positive roots.   
%   The $(\bg,K)$-module $W(\lambda)$  is cohomologically induced from a one dimensional representation of $\bq =\bl \oplus \bn $
%which is  trivial on $SO_o(p,q-2)$ and a character with parameter $\lambda $ on the torus SO(2).  Here $\lambda $ is a positive integer . {\color{magenta} Correct !!!! {\cite
%S} ? see [K] page 27)}. We denote this $(\bg,K)$-module by $W(\lambda)$.  Since we assume that $p>3$  L is noncompact  and thus the $(\bg,K)$-module $W(\lambda)  $ is not tempered.

%\medskip
 %Following the notation of Y.~Oshima [O] and  with the notation of [KS3] 

%*****************************************************

 \section{Another perspective and some speculations. }
  In this section we rephrase our previous results and reexamine them in the 
context of Arthur packets and Arthur Vogan packets. 
  
 To simplify the notation and the considerations we assume  in this section that the groups 
$G$, $G^1$ and $G^2$ have discrete series representations, i.e. that p and q are even and that 
$p \leq q$. Furthermore we assume as in the earlier sections that 
  $ 4 \leq p,$

 \bigskip
 \subsection{Generalities}
 An Arthur Vogan parameter is a  a homomorrphism of the Weil group in the Langlands dual group of G which commutes with a  homorphism of $SL(2,\bC)$. This defines a unipotent orbit and  a one dimensional representation of the Levi subgroup of a parabolic subgroup of the complexification $G_\bC$ of G. We obtain a  family of one dimensional representations  of $\theta$-stable parabolic subgroups of a family of pure inner forms  of $G_\bC$. By cohomological induction we obtain a family of irreducible representations of pure inner forms compatible G.  For details, definitions and results  see the article by  J.~Adams and J.~Johnson \cite{A-J} or the book by Adams-Barbasch-Vogan \cite{A-B-V}.
 
 \medskip
 The family of pure inner forms of  $SO_o(p,q) $ consists of the groups $SO_o(p-2r, q+2r),$  $q<2r <p$ and $SO_o(p+2s,q-2s)$. Similarly we obtain a family  of pure inner forms compatible with $G^1$ and $G^2$ and observe that these inner classes are disjoint.
 
 \medskip
 
{\bf  Examples:} \\
 If both $p$ and $q$ are even then the pure inner class of pure inner forms of 
$SO_o(p,q)$ contains both the compact groups $SO(p+q,0)$ and $SO(0,p+q).$\\
  The class of pure inner forms of $SO_0(2n,2m+1) $ contains the compact group $SO(0,2n+2m +1)$ and 
and   $SO_o(2n+2m,1)$ the orthogonal group of real rank one.

 \medskip
 
   Using lemma 2.10 of J.~Adams and J.~Johnson \cite{A-J}  we can show, that  if $G=SO_o(p,q) $ with both 
$p,q$ even and $4\leq p \leq q$ then the number of inequivalent irreducible representations in an  
Arthur packet containing ${W}(\lambda)^\infty$  is 2 as follows:
   We denote by $T$ a compact Cartan subgroup of $SO(p,q)$ and $\bt_\bC$ its complexifed Lie algebra. 
The Lie algebra of $G_\bC$ is noted by  by $\bg_\bC$ and  
$\bl_\bC = \bl \otimes \bC = \bC \oplus so(p-2,q)$. The corresponding Weyl groups are $W(G,T)$, $W(\bg_\bC, \bt_\bC )$ and $W(\bl_\bC,\bt_\bC)$. The representations in an Arthur packet are parametrized by the double cosets
   \[ {\mathbb S}=W(G,T) \backslash W(\bg_\bC) /W(\bl_\bC,\bt_\bC) \]
   
   \begin{lemma}
   $\mathbb S$ has 2 disjoint double cosets.
   \end{lemma}
   \proof
 Note that the coset $  W(\bg_\bC) /W(\bl_\bC,\bt_\bC) $ has $ m=2[\frac{p+q}{2}]$ elements and has representants $P_i$ acting on a basis $e_i$ of $\bC^m$ by permuting  $e_1$ and $e_i$ and by sign changes composed with permutations. W(G,T) is the Weyl group of $SO(p) \times SO(q)$. It has 2 orbits  on $  W(\bg_\bC) /W(\bl_\bC,\bt_\bC) $ with representants Identity and $P_1$.
  \qed
  
  \medskip
  \begin{remark}
  A similar proof shows that $\mathbb S$ contains exactly 2 disjoint double cosets if $p$ or $q$ is even and 
   $4\leq p \leq q$.
  \end{remark}
  
  {\bf Examples:}\\
  For $q=1, p=2m, \ m>I $ the results in \cite{K-S2} show that in this case the Arthur packet containing $W(\lambda)^\infty$ contains only one element, a representation of height one.\\
  If $p=0$, then the Arthur packet contains exactly one finite dimensional spherical 
subrepresentation of $L^2(SO(0,q)/SO(0,q-1))$.

  \medskip
   Consider now a maximal parabolic subgroup $P_\bC$ of $SO(p+q,\bC)$ with a Levi $L_\bC=SO(p+q-2) ~\bC^*$.  If  $4\leq p \leq q$ it defines two  $\theta$-stable parabolic subgroups $P_1$ and $P_2$ of $SO_o(p,q)$ with Levi 
 $L_1=SO(p,q-2) SO(0,2)$ respectively $L_2=SO(2,0)SO(p-2,q)$, which are both in the same family of pure inner forms. 
   If two irreducible representations  are cohomological induced from characters of $P_1$ and $P_2$ and have the same infinitesimal character  then they are in the same Vogan packet \cite{A-J}.

   \medskip
 We  denote the representation in the Arthur packet corresponding to the double coset of the identity  is denoted by  ${W}(\lambda)^\infty$  and the representation corresponding to the other coset by ${W}_{as}(\lambda)^\infty$. ${W}(\lambda)^\infty$ is in the discrete spectrum of $L^2(SO_o(p+q)/SO_o(p,q-1)$ whereas the representation ${W}_{s}(\lambda)^\infty$ is  also a Flensted Jensen representation in the discrete spectrum of  $L^2(SO_o(p,q) /SO_o(p-1,q))$.
 
 \medskip
 An Arthur Vogan packet containg the Arthur packet $\{ {W}_{as}(\lambda)^\infty$. $\hat{W}(\lambda) ^\infty \}$
 is a disjoint union of the Arthur packets of the Flensted-Jensen representations of 
the pure inner forms of $SO(p,q)$.
  %If they have the same infinitesimal character then they are in the same 
 %Arthur packet.
 \medskip
 
 {\bf Example} The case $G=SO_0(3,3)$ was partially discussed in \cite{O-S}. We consider there a 
representation cohomological induced from the trivial character of $L=SO(2,0)SO_o(1,3)$, i.e. a 
Flensted-Jensen representation of $SO_0(3,3)/SO_o(2,3)$, i.e. in our notation $W(0)^\infty_{as}$. 
The other representation is the Arthur packet is cohomological induced from the trivial representation of 
$L=SO(0,2)SO_o(3,1)$.\\
 The pure inner forms of $SO_0(3,3)$ are $SO_o(5,1)$ and $SO_o(1,5)$. We are considering the representations in the discrete spectrum of $SO_0(1,5)/SO_o(1,4)$ respectively $SO_o(5,1)/SO(4,1)$, since a symmetric space has no discrete spectrum.  These representations are cohomological induced from $L=SO_o(1,3)SO(0,2)$ and $L=SO(2,0)SO_0(3,1)$. \\
 Thus the Arthur packet consists of  2  representations and the Arthur Vogan packet contains  4 representations.
 
 \bigskip
\subsection{Restriction of representations in Arthur packets} We fix a regular infinitesimal character
 \[(\lambda,0,\dots ,0) + \rho \]
of $SO(p+q,\bC)$  satisfying conditions 8.8 and 8.9 in \cite{S2}. 
% The Arthur packet  containing   $\hat{W}(\lambda)$ in the discrete spectrum of  contains also another  representation $\hat{W}_{as}(\lambda)$ cohomological induced from the other $\theta $ stable parabolic subgroup. Since the parabolic subgroups are not conjugate under K these 2 representations are not isomorphic.
%(I believe looking at Schlichtkrulls formulas for the Langlands parameter) that the second representation is a anti de sitter rep. Call them Pi_s and \Pi_a.  

%We consider the restriction of the representations  $\hat{W}(\lambda)$ and $\hat{W}_{as}(\lambda)$ to the subgroups $G^1= SO(p-1,q)$ and $G^2 = SO(p,q-1) $.

 $L$-functions and  trace formula considerations strongly suggest  that this symmetry breaking  depends  on the inner class of the groups $G^i$  and the Arthur packet of representations.
 Our computations and the results of \cite{K-O}  support this. Recall that the restriction of a representations to a subgroup H is called admissible if it is direct sum of representations of H.
 
 \begin{prop} 
 
 \begin{enumerate} \item The restriction to $G^1$ of
 exactly one of the representations   $\{ {W}(\lambda)^\infty, \ {W}_{as}(\lambda)^\infty \} $ in the Arthur packet is admissible, the other one is not admissible.
\item  The restriction to $G^2$ of
 exactly one of the representations   $\{ {W}(\lambda)^\infty, \ {W}_{as}(\lambda)^\infty \} $ in the Arthur packet is admissible, the other one is not admissible.
 \item Each representation in the Arthur packet is admissible for exactly one of the groups $G^i$, $ i\in \{ 1, \ 2 \}$
 
 \end{enumerate}
\end{prop}
\proof
For  ${W}(\lambda)^\infty$ this follows from the branching theorem and for ${W}_{as}(\lambda)^\infty$ it follows from the table in \cite{K-O}.  More precise information about the restriction  of ${W}_{as}(\lambda)^\infty$ in the unitary category is contained  in \cite{Ko93} supporting the conjecture.
\qed

\medskip
\begin{remark}
The complete branching law for the restriction of $W_{as}(\lambda)$ to $G^2$ can be found in \cite{Ko93} theorem 3.3 and information about the discrete spectrum of the restriction of the unitary representation  with Harish-Chandra-module $W_{as}(\lambda )$ to $G^1$  are in \cite{Ko02}.
\end{remark}

\bigskip
\subsection{} Based on our computations we expect a stronger result similar to the
 famous Gross-Prasad conjecture. Our assumptions imply that for the groups $G^i$, i=1,2 we have an 
Arthur packets with 2 representations $ \{ {V}(\nu)^\infty, {V}_{as}(\nu)^\infty \}$ , 
respectively  $\{ U(\mu)^\infty, U_{as}(\mu)^\infty \}$.

{\bf Conjecture for Arthur packets} \\
Suppose  that $\pi_1$ is a representations of $G^1$ in an Arthur  packet  \[{\mathcal C_{p-1,q}} =
\{ {V}(\nu)^\infty, {V}_{as}(\nu)^\infty \} .\]
\begin{enumerate}
\item If 
                  \[    \mbox{Hom} _{G^1} ( W(\lambda)^\infty _{|G^1}, \pi_1^\infty ) \not = 0 \]
then 
                  \[     \mbox{Hom} _{G^1} ( {W}_{as}(\lambda)^\infty _{| G^1}, \pi_1^\infty ) = 0  \]

\item if        \[      \mbox{Hom} _{G^1} ( {W}_{as}(\lambda)^\infty _{|G^1}, \pi_1^\infty ) \not = 0  \]

then        \[         \mbox{Hom} _{G^1} ( {W}(\lambda)^\infty _{G^1}, \pi_1^\infty ) = 0  \]
\end{enumerate}

There is an analogous statement for the group  $G^2$.

\medskip
{\bf Definition}
If ${G}$ is a pure inner form of $SO_o(p,q)$ and ${G}^i $ is a pure inner form  of $G^i$ we call (${G}^i$  ${G}$) a relevant pair if ${G}^i$ is a subgroup of ${G}$. 

\medskip
{\bf Example:} The pure inner forms of $SO_o(3,3)$ are
\[  SO_0(1,5)  \ \ \ \ \ \ \ SO_0(3,3) \ \ \ \ \ \ \ \ SO_o(5,1) \]
The pure inner forms of $G^1$ are 
\[ SO_o(0,5) \ \ \ \ \ \ \ \  SO_o(2,3) \ \ \ \ \ \ \ SO_o(4,1)\]
and relevant pairs
\[ (SO_o(1,5),SO(5))  \ \ \ \ \ (SO_o(3,3),SO_o(2,3) ) \ \ \ \ \ (SO_o(5,1),SO_0(4,1)) \]

\medskip

{\bf Conjecture for Arthur-Vogan packets} \\
Let ${\mathcal AV} _G$ be the Arthur Vogan packet  which contains ${W}(\lambda)^\infty$ and ${\mathcal AV}_{G^1}$ the Arthur packet which contains ${V}(\nu)^\infty$. Suppose that $\underline{G}$ is a pure inner form of $G=SO_o(p,q)$ and $\underline{G}^1 $ is a pure inner form of $G^1=SO_o(p-1,q)$ and that $\underline{G}$ and $\underline{G}^1$ is a relevant pair.

 Assume that  $W$ and $W_a $ are inequivalent representations of $\underline{G}$ in the Arthur packet $(\mathcal AV)_{{G}}$ and $V$ a representation of
 $\underline{G}^1$ in the Arthur-Vogan packet  $ (\mathcal AV)_{G^1}$.
 
 \begin{enumerate}
 \item If 
                  \[    \mbox{Hom} _{\underline{G}^1} (W^\infty,V^\infty) ) \not = 0 \]
then 
                  \[     \mbox{Hom} _{\underline{G}^1} ( W_a^\infty, V^\infty ) = 0  \]

\item If        \[      \mbox{Hom} _{\underline{G}^1} (W_a^\infty , V^\infty ) \not = 0  \]

then        \[         \mbox{Hom} _{\underline{G}^1} ( W^\infty ,V^\infty ) = 0  \]
\end{enumerate}

We expect that  an analogous statement holds for  ${\mathcal AV}_{G^2}$.

\begin{remark}
The conjecture for Arthur-Vogan packets is similar but weaker than  the Gross Prasad conjectures for Vogan packets of  tempered representations \cite{GP},  \cite{A1}.

\end{remark}

\begin{remark}
If $G=SO_o(n+1,1) $ then  the Arthur packet      has only one representation,  $G^1=SO(n+1)$ is compact and $G^2  = SO_o(n,1) $ and the results follow from \cite{K-S2}.             
\end{remark}

% All the representations of G,G' which we considered are in the same Arthur Vogan packet. So we conclude
 
% \begin{theorem}
% Multiplicity one fails for Arthur Vogan packet, which are not also Vogan -packets. 
% \end{theorem}

% \newpage
% \section{An example : G= O(3,3),  L=SO(2)O(1,3),  $G^1$= O(2,3), $G^2$ = O(3,2)}

\ 
%\newpage

\end{document}